\documentclass[12pt,fleqn]{article}
\usepackage{amsfonts,amsmath,amsthm}

\input xy
\xyoption{all}
\CompileMatrices

\newcommand{\bdry}[1]{\partial #1}
\newcommand{\bigset}[1]{\big\{#1\big\}}
\newcommand{\C}[1]{C_{#1}}
\newcommand{\Cl}[1]{C_{#1 1}}
\newcommand{\Cu}[1]{C_{#1 2}}
\newcommand{\dhalf}{\displaystyle \frac{1}{2}}
\newcommand{\dist}[2]{\text{dist}(#1,#2)}
\newcommand{\etabar}{\overline{\eta}}
\newcommand{\FUCIK}{{F}u\v c\'\i k }
\newcommand{\Fucik}{Fu\v c\'\i k }
\newcommand{\halfquad}{\hspace{0.08in}}
\newcommand{\halfthinspace}{\hspace{0.016in}}
\newcommand{\ip}[2]{\left(#1,#2\right)}
\newcommand{\norm}[2]{\|#1\|_{#2}}
\newcommand{\R}{\mathbb R}

\newcommand{\seq}[1]{\left(#1\right)}
\newcommand{\set}[1]{\left\{#1\right\}}
\newcommand{\Z}{\mathbb Z}

\DeclareMathOperator{\divg}{div}
\DeclareMathOperator{\sgn}{sgn}

\newenvironment{enumeratepr}{\begin{enumerate}

}{\end{enumerate}}
\newenvironment{enumeratepre}{\begin{enumerate}

}{\end{enumerate}}

\newtheorem{corollary}{Corollary}[section]
\newtheorem{lemma}[corollary]{Lemma}
\newtheorem{proposition}[corollary]{Proposition}
\newtheorem{theorem}[corollary]{Theorem}

\theoremstyle{remark}

\numberwithin{equation}{section}

\title{\vspace{-0.5in}
       {\normalsize To appear in J. Math. Anal. Appl.}\\
       \vspace{0.2in}
       \LARGE \bf Some Remarks on the \Fucik Spectrum of the
$p$-Laplacian and Critical Groups}
\author{\large \bf Norman Dancer\\ \large School of Mathematics
and Statistics\\ \large University of Sydney\\ \large NSW 2006,
Australia\\[0.05in] {\em E-mail:} normd@maths.usyd.edu.au\\ \\
\large \bf Kanishka Perera\thanks{The second author gratefully
acknowledges the support and hospitality of the School of
Mathematics and Statistics at the University of Sydney where
this work was done.}\\ \large Department of Mathematical
Sciences\\ \large Florida Institute of Technology\\ \large
Melbourne, FL 32901, U. S. A.\\[0.05in]
{\em E-mail:} kperera@winnie.fit.edu}
\date{}

\begin{document}

\maketitle

\section{Introduction} \label{S1}

Let $\Omega$ be a bounded domain in $\R^n,\, n \ge 1$. For $1 <
p < \infty$, the \Fucik spectrum of the $p$-Laplacian on
$W^{1,\, p}_0(\Omega)$ is defined as the set $\Sigma_p$ of those
points $(a,b) \in \R^2$ such that
\begin{equation} \label{1}
- \Delta_p\, u = a\, (u^+)^{p-1} - b\, (u^-)^{p-1}, \quad u \in
W^{1,\, p}_0(\Omega)
\end{equation}
has a nontrivial solution. Here $\Delta_p\, u = \divg \big(|\nabla
u|^{p-2}\, \nabla u\big)$ and $u^\pm = \linebreak \max \set{\pm\,
u, 0}$. It is known that the first eigenvalue $\lambda_1$ of $-
\Delta_p$ on $W^{1,\, p}_0(\Omega)$ is positive, simple, and
admits a positive eigenfunction $\varphi_1 \in W^{1,\,
p}_0(\Omega) \cap C^1(\Omega)$ (see Lindqvist \cite{Lin1}), so
$\Sigma_p$ contains the two lines $\lambda_1 \times \R$ and $\R
\times \lambda_1$. A first nontrivial curve $\C{2}$ in $\Sigma_p$
through $(\lambda_2,\lambda_2)$, where $\lambda_2$ is the second
eigenvalue of $- \Delta_p$, was recently constructed and
variationally characterized by a mountain-pass procedure in
Cuesta, de Figueiredo, and Gossez \cite{CudeFiGo}. It was shown
there that $\C{2}$ is continuous, strictly decreasing, and
asymptotic to $\lambda_1 \times \R$ and $\R \times \lambda_1$ at
infinity.

As is well-known, solutions of \eqref{1} are the critical
points of the $C^1$ functional
\begin{equation} \label{2}
I_{(a,b)}(u) = \int_\Omega |\nabla u|^p - a\, (u^+)^p - b\,
(u^-)^p, \quad u \in X = W^{1,\, p}_0(\Omega).
\end{equation}
When $(a,b) \not \in \Sigma_p$, the origin is an isolated critical
point of $I_{(a,b)}$ and hence the critical groups
$C_\ast(I_{(a,b)},0)$ are defined. Some of these groups were
computed in Dancer \cite{Da2} and Perera and Schechter
\cite{PeSc9, PeSc2} in the linear case $p = 2$. The purpose of
the present paper is to extend some of these computations to the
quasilinear case $p \ne 2$. Let $\Cl{1} = \bigl((-
\infty,\lambda_1] \times \lambda_1\bigr) \cup \bigl(\lambda_1
\times (- \infty,\lambda_1)\bigr)$ and $\Cu{1} = \bigl(\lambda_1
\times [\lambda_1,+ \infty)\bigr) \cup \bigl((\lambda_1,+ \infty)
\times \lambda_1\bigr)$. We shall show that

\begin{theorem} \label{T2}
\begin{enumeratepre}
\item \label{10} If $(a,b)$ lies below $\Cl{1}$, then $C_q(I_{(a,b)},0)
= \delta_{q0}\, \Z$.
\item \label{11} If $(a,b)$ lies between $\Cl{1}$ and $\Cu{1}$, then
$C_q(I_{(a,b)},0) = 0 \halfquad \forall q$.
\item \label{12} If $(a,b)$ lies between $\Cu{1}$ and $\C{2}$, then
$C_q(I_{(a,b)},0) = \delta_{q1}\, \Z$.
\item \label{13} If $(a,b) \notin \Sigma_p$ lies above $\C{2}$, then
$C_q(I_{(a,b)},0) = 0,\; q = 0,1$.
\end{enumeratepre}
\end{theorem}

Denote by $\widetilde{J}_s$ the restriction of the functional
\begin{equation} \label{3}
J_s(u) = \int_\Omega |\nabla u|^p - s\, (u^+)^p, \quad u \in X
\end{equation}
to the $C^1$ manifold
\begin{equation} \label{4}
S = \set{u \in X : \int_\Omega |u|^p = 1}.
\end{equation}
As noted in \cite{CudeFiGo}, the points in $\Sigma_p$ on the line
parallel to the diagonal $a = b$ and passing through $(s,0)$ are
exactly of the form $(s + c,c)$ with $c$ a critical value of
$\widetilde{J}_s$. As we will see in Section \ref{S2}, the
critical groups of $I_{(a,b)}$ are related to the homology groups
of the sublevel sets of $\widetilde{J}_{a-b}$. Note that the
standard second deformation lemma cannot be used in determining
the structure of these sublevel sets as the manifold $S$ is not
of class $C^{1,1}$ when $p < 2$. We will overcome this difficulty
by using a deformation lemma on a $C^1$ manifold from Ghoussoub
\cite{Gh2} and the Ekeland's variational principle to show that
$H_\ast(\widetilde{J}_s^\alpha) \cong
H_\ast(\widetilde{J}_s^\beta)$ if $\widetilde{J}_s$ has no
critical values in $[\alpha,\beta]$. We will also use this to
prove the following homotopy invariance result for
$C_\ast(I_{(a,b)},0)$, which does not follow from the standard
homotopy invariance theorem for critical groups again because of
our very limited smoothness.

\begin{proposition} \label{T6}
If $(a_0,b_0)$ and $(a_1,b_1)$ can be joined by a curve that does
not intersect $\Sigma_p$, then
\begin{equation} \label{4.7}
C_\ast(I_{(a_0,b_0)},0) \cong C_\ast(I_{(a_1,b_1)},0).
\end{equation}
\end{proposition}

In the process of proving Theorem \ref{T2}, we will also show that
the curve $\C{2}$ has the following topological property. For
$b > \max \set{\lambda_1 - s, \lambda_1}$, the set
\begin{equation} \label{4.701}
{\cal O}_b = \set{u \in S : \widetilde{J}_s(u) < b}
\end{equation}
is path-connected if and only if the point $(s + b,b)$ lies above
$\C{2}$. This was proved in Dancer and Du \cite{DaDu} for $p = 2$.

As an application, we consider the quasilinear elliptic boundary
value problem
\begin{equation} \label{39}
- \Delta_p\, u = f(x,u), \quad u \in X
\end{equation}
where $f$ is a Carath\'{e}odory function on $\Omega \times \R$
such that
\begin{equation} \label{40}
f(x,t) = \begin{cases}
a_0\, (t^+)^{p-1} - b_0\, (t^-)^{p-1} + o(|t|^{p-1}) & \text{as }
t \to 0,\\[5pt]
a\, (t^+)^{p-1} - b\, (t^-)^{p-1} + o(|t|^{p-1}) & \text{as } |t|
\to \infty
\end{cases}
\end{equation}
uniformly in $x$. We shall prove

\begin{theorem} \label{T9}
\begin{enumeratepre}
\item \label{60} If the points $(a_0,b_0),\, (a,b)$ lie on opposite
sides of \linebreak $\lambda_1 \times \R$, then \eqref{39} has a
positive solution.
\item \label{61} If $(a_0,b_0),\, (a,b)$ lie on opposite sides of
$\R \times \lambda_1$, then \eqref{39} has a negative solution.
\item \label{59} If $(a_0,b_0),\, (a,b) \notin \Sigma_p$ lie on
opposite sides of $\C{2}$, then \eqref{39} has a nontrivial
solution.
\item \label{62} If $(a_0,b_0),\, (a,b)$ lie on opposite sides of
either $\Cl{1}$ or $\Cu{1}$, then there is a fixed-sign solution,
and if they lie on opposite sides of both $\Cl{1}$ and $\Cu{1}$,
then there are a positive solution and a negative solution.
\item \label{63} If $(a_0,b_0) \notin \Sigma_p$ is above $\C{2}$
and $(a,b)$ is below $\Cl{1}$, then there is a third nontrivial
solution.
\end{enumeratepre}
\end{theorem}

In particular, \eqref{39} has a nontrivial solution if $a_0 =
b_0,\, a = b$ are not in $\sigma(- \Delta_p)$ and $\lambda_1$ or
$\lambda_2$ lies between $a_0$ and $a$, which is a special case of
the well-known Amann-Zehnder theorem when $p = 2$ (see
\cite{AmZe}).

Solutions of \eqref{39} are the critical points of
\begin{equation} \label{41}
\Phi(u) = \int_\Omega |\nabla u|^p - p\, F(x,u), \quad u \in X
\end{equation}
where $F(x,t) = {\displaystyle \int_0^t} f(x,s)\, ds$.
In the nonresonance case $(a_0,b_0),\, (a,b) \notin \Sigma_p$,
one would expect that
\begin{equation} \label{42}
C_\ast(\Phi,0) \cong C_\ast(I_{(a_0,b_0)},0), \qquad
C_\ast(\Phi,\infty) \cong C_\ast(I_{(a,b)},0).
\end{equation}
This, however, seems difficult to prove since $\Phi$ is not
$C^{1,1}$ in general. To get around this difficulty we will
construct a perturbed functional $\widetilde{\Phi}$ that has the
same critical points as $\Phi$, with $\widetilde{\Phi}(u) =
I_{(a_0,b_0)}(u)$ for $\norm{u}{}$ small and $\widetilde{\Phi}(u)
= I_{(a,b)}(u)$ for $\norm{u}{}$ large. Then \eqref{42} would
hold with $\Phi$ replaced by $\widetilde{\Phi}$. Proof of Theorem
\ref{T9} is given in Section \ref{S3}.

\section{Critical Group Computations} \label{S2}

We work with the constrained functional $\widetilde{J}_s$. Since
$\Sigma_p$ is clearly symmetric with respect to the diagonal, we
may assume that $s \ge 0$. It was shown in Cuesta, de Figueiredo,
and Gossez \cite{CudeFiGo} that $\widetilde{J}_s$ satisfies (PS),
\begin{enumeratepr}
\item \label{4.3} $\varphi_1$ is a global minimum of
$\widetilde{J}_s$ with $\widetilde{J}_s(\varphi_1) = \lambda_1 -
s$, and the corresponding point $(\lambda_1,\lambda_1 - s)$ in
$\Sigma_p$ lies on the vertical line through
$(\lambda_1,\lambda_1)$,
\item \label{4.4} $- \varphi_1$ is a strict local minimum of
$\widetilde{J}_s$ with $\widetilde{J}_s(- \varphi_1) = \lambda_1$,
and the corresponding point $(s + \lambda_1,\lambda_1)$ in
$\Sigma_p$ lies on the horizontal line through
$(\lambda_1,\lambda_1)$,
\item \label{4.5} denoting by $\Gamma$ the family of all continuous
paths in $S$ joining $\varphi_1$ and $- \varphi_1$,
\begin{equation} \label{4.2}
c(s) := \inf_{\gamma \in \Gamma} \max_{u \in \gamma([0,1])}
\widetilde{J}_s(u)
\end{equation}
is the first critical value of $\widetilde{J}_s$ that is $>
\lambda_1$, so $\C{2} = \{(s + c(s),c(s)) : s \ge 0\} \cup
\{(c(s),s + c(s)) : s \ge 0\}$ is the first nontrivial curve in
$\Sigma_p$.
\end{enumeratepr}

\begin{lemma} \label{T4}
Denoting by $K_c$ the set of critical points of $\widetilde{J}_s$
with critical value $c$, for $s > 0$, $K_{\lambda_1 - s} =
\set{\varphi_1},\, K_{\lambda_1} = \set{- \varphi_1}$, and
$\widetilde{J}_s$ has no critical values in $(\lambda_1 - s,\lambda_1)$.
\end{lemma}

\begin{proof}
If $u \in K_{\lambda_1}$, then
\begin{equation} \label{9.3}
- \Delta_p\, u = (\lambda_1 + s)\, (u^+)^{p-1} - \lambda_1\,
(u^-)^{p-1}.
\end{equation}
Since any eigenfunction of $- \Delta_p$ associated with an
eigenvalue $> \lambda_1$ changes sign, $u < 0$ somewhere. Let
$\Omega^- = \set{x \in \Omega : u(x) < 0}$. Then
\begin{equation} \label{9.4}
- \Delta_p\, u^- = \lambda_1\, (u^-)^{p-1} \quad \text{in }
\Omega^-
\end{equation}
and $u^- \in W^{1,\, p}_0(\Omega^-)$ by Lemma 5.6 of
\cite{CudeFiGo}. Since $u^- > 0$ in $\Omega^-$, $\lambda_1$ must
also be the first Dirichlet eigenvalue of $- \Delta_p$ on
$\Omega^-$, so it follows from Lemma 5.7 of \cite{CudeFiGo} that
$\Omega^- = \Omega$. Thus $u = - \varphi_1$.

We complete the proof by showing that if $u \in K_c$ with $c <
\lambda_1$, then $c = \lambda_1 - s$ and $u = \varphi_1$. We have
\begin{equation} \label{9.5}
- \Delta_p\, u = (s + c)\, (u^+)^{p-1} - c\, (u^-)^{p-1}.
\end{equation}
If $\Omega^- = \set{x \in \Omega : u(x) < 0} \neq \emptyset$, then
$u^- \in W^{1,\, p}_0(\Omega^-)$ satisfies
\begin{equation} \label{9.6}
- \Delta_p\, u^- = c\, (u^-)^{p-1} \quad \text{in } \Omega^-,
\end{equation}
so $\lambda_1(\Omega^-) = c < \lambda_1$, contradicting Lemma 5.7
of \cite{CudeFiGo}. Thus $u \ge 0$ and satisfies
\begin{equation} \label{9.7}
- \Delta_p\, u = (s + c)\, u^{p-1}. \qed
\end{equation}
\renewcommand{\qed}{}
\end{proof}

Note that
\begin{equation} \label{4.1}
I_{(a,b)}|_S = \widetilde{J}_{a-b} - b,
\end{equation}
so the sublevel sets
\begin{equation} \label{4.11}
I_{(a,b)}^\alpha = \Big\{u \in X : I_{(a,b)}(u) \le \alpha\Big\},
\qquad \widetilde{J}_s^\alpha = \Big\{u \in S : \widetilde{J}_s(u)
\le \alpha\Big\}
\end{equation}
are related by
\begin{equation} \label{4.13}
I_{(a,b)}^\alpha \cap S = \widetilde{J}_{a-b}^{\alpha + b}.
\end{equation}

\begin{lemma} \label{T1}
If $(a,b) \notin \Sigma_p$, then
\begin{equation} \label{7}
C_q(I_{(a,b)},0) \cong \begin{cases}
\delta_{q0}\, \Z & \text{if } \widetilde{J}_{a-b}^b =
\emptyset,\\[5pt] \widetilde{H}_{q-1}(\widetilde{J}_{a-b}^b) &
\text{otherwise},
\end{cases}
\end{equation}
where $\widetilde{H}_\ast$ denote reduced homology groups.
\eqref{7} also holds with $\widetilde{J}_{a-b}^b$ replaced by
${\cal O}_b = \set{u \in S : \widetilde{J}_{a-b}(u) < b}$.
\end{lemma}

\begin{proof}
Since $0$ is the only critical point of $I_{(a,b)}$,
\begin{equation} \label{8}
C_q(I_{(a,b)},0) = H_q(I_{(a,b)}^0,I_{(a,b)}^0 \setminus \set{0}).
\end{equation}
As $I_{(a,b)}$ is positive homogeneous, $I_{(a,b)}^0$ is
contractible and $I_{(a,b)}^0 \setminus \set{0}$ is homotopic to
$I_{(a,b)}^0 \cap S$, so \eqref{7} follows from the reduced
homology sequence of the pair $(I_{(a,b)}^0,I_{(a,b)}^0 \setminus
\set{0})$ and \eqref{4.13}. The second statement is proved
similarly.
\end{proof}

We will say that a closed subset $\cal O$ of
$\widetilde{J}_s^\beta$ is isolated if $\cal O$ and
$\widetilde{J}_s^\beta \setminus {\cal O}$ cannot be connected by
a path in $\widetilde{J}_s^\beta$. Since $S$ is not $C^{1,1}$ when
$p < 2$, our next lemma does not follow from the standard second
deformation lemma.

\begin{lemma} \label{T5}
If ${\cal O} \subset \widetilde{J}_s^\beta$ is isolated and
$\widetilde{J}_s$ has no critical points in \linebreak
$\widetilde{J}_s^{-1}([\alpha,\beta]) \cap {\cal O}$, then
\begin{equation} \label{4.14}
H_\ast({\cal O},\widetilde{J}_s^\alpha \cap {\cal O}) = 0.
\end{equation}
In particular, if $\widetilde{J}_s$ has no critical values in
$[\alpha,\beta]$, then
\begin{equation} \label{4.6}
H_\ast(\widetilde{J}_s^\beta,\widetilde{J}_s^\alpha) = 0.
\end{equation}
\end{lemma}

Before the proof, two corollaries. The first should be compared
with Lemma 3.6 of \cite{CudeFiGo}.

\begin{corollary} \label{T7}
If $\cal O$ is a nonempty isolated subset of
$\widetilde{J}_s^\beta$, then $\cal O$ contains a critical point
of $\widetilde{J}_s$.
\end{corollary}

\begin{proof}
Taking $\alpha < \inf \widetilde{J}_s({\cal O})$ gives
\begin{equation}
H_0({\cal O},\widetilde{J}_s^\alpha \cap {\cal O}) = H_0({\cal O})
\ne 0. \qed
\end{equation}
\renewcommand{\qed}{}
\end{proof}

\begin{corollary} \label{T8}
If ${\cal O} \subset \widetilde{J}_s^\beta$ is isolated and the
only critical point of $\widetilde{J}_s$ in $\cal O$ is a local
minimizer $u_0$ such that
\begin{equation} \label{4.71}
\widetilde{J}_s(u) > \widetilde{J}_s(u_0) \quad \forall u \in
{\cal O} \setminus \set{u_0}, \qquad \inf_{u \in
\bdry{\widetilde{B}_\varepsilon(u_0)}} \widetilde{J}_s(u) >
\widetilde{J}_s(u_0)
\end{equation}
for some $\varepsilon > 0$ such that
$\widetilde{B}_\varepsilon(u_0) = \set{u \in S : \norm{u - u_0}{}
\le \varepsilon} \subset {\cal O}$, then
\begin{equation} \label{4.73}
\widetilde{H}_\ast({\cal O}) = 0.
\end{equation}
\end{corollary}

\begin{proof}
Take $\widetilde{J}_s(u_0) < \alpha < \inf
\widetilde{J}_s(\bdry{\widetilde{B}_\varepsilon(u_0)})$. Then
${\cal O}' = \widetilde{J}_s^\alpha \cap {\cal O} \setminus
\widetilde{B}_\varepsilon(u_0)$ is isolated in
$\widetilde{J}_s^\alpha$ and $\widetilde{J}_s$ has no critical
points in ${\cal O}'$, so ${\cal O}' = \emptyset$ by Corollary
\ref{T7}, i.e., $\widetilde{J}_s^\alpha \cap {\cal O} \subset
\widetilde{B}_\varepsilon(u_0)$. Thus we have the commutative
diagram
\begin{equation} \label{23}
\xymatrix{{\widetilde{H}}_\ast(\widetilde{J}_s^\alpha \cap {\cal
O}) \ar[r] \ar[dr]_{i_\ast} &
{\widetilde{H}}_\ast(\widetilde{B}_\varepsilon(u_0)) \ar[d]\\
& {\widetilde{H}}_\ast({\cal O})}
\end{equation}
induced by inclusions, where $i_\ast$ is an isomorphism by Lemma
\ref{T5} as $\widetilde{J}_s$ has no critical points in
$\widetilde{J}_s^{-1}([\alpha,\beta]) \cap {\cal O}$. Since
$\widetilde{H}_\ast(\widetilde{B}_\varepsilon(u_0)) = 0$, the
conclusion follows.
\end{proof}

\begin{proof}[Proof of Lemma \ref{T5}]
If $[z] \in H_q({\cal O},\widetilde{J}_s^\alpha \cap {\cal O})$ is
nontrivial, set
\begin{equation} \label{4.8}
c := \inf_{z' \in [z]} \max_{u \in |z'|} \widetilde{J}_s(u)
\end{equation}
where $|z'|$ denotes the support of the singular $q$-chain $z'$.
Clearly, $c \in [\alpha,\beta]$. We will show that
$\widetilde{J}_s$ has a critical point in $\widetilde{J}_s^{-1}(c)
\cap {\cal O}$.

We follow the proof of Theorem 3.2 in Ghoussoub \cite{Gh2}.
Consider the subspace $\cal L$ of $C([0,1] \times S; S)$
consisting of all continuous deformations $\eta$ such that
\begin{enumeratepr}
\item \label{4.81} $\eta(0,u) = u \quad \forall u \in S$,
\item \label{4.82} $\sup \set{\rho(\eta(t,u),u) : (t,u) \in [0,1] \times S}
< + \infty$ where $\rho$ is the Finsler metric on $S$,
\item \label{4.83} $\widetilde{J}_s(\eta(t,u)) \le \widetilde{J}_s(u) \quad
\forall (t,u) \in [0,1] \times S$.
\end{enumeratepr}
$\cal L$ equipped with the metric
\begin{equation} \label{4.9}
\delta(\eta,\eta') = \sup \set{\rho(\eta(t,u),\eta'(t,u)) : (t,u)
\in [0,1] \times S}
\end{equation}
is a complete metric space. For any $\eta \in {\cal L}$, the
assumption that there is no path in $\widetilde{J}_s^\beta$
joining $\cal O$ and $\widetilde{J}_s^\beta \setminus {\cal O}$
together with \eqref{4.83} above imply that, for each $t \in
[0,1]$, the restriction of $\eta(t,\cdot)$ to $\cal O$ is a map of
the pair $({\cal O},\widetilde{J}_s^\alpha \cap {\cal O})$ (i.e.,
$\eta(t,{\cal O}) \subset {\cal O}$ and
$\eta(t,\widetilde{J}_s^\alpha \cap {\cal O}) \subset
\widetilde{J}_s^\alpha \cap {\cal O}$). Thus, \eqref{4.81} implies
that $\eta(1,\cdot) : ({\cal O},\widetilde{J}_s^\alpha \cap {\cal
O}) \to ({\cal O},\widetilde{J}_s^\alpha \cap {\cal O})$ is
homotopic to the identity on $({\cal O},\widetilde{J}_s^\alpha
\cap {\cal O})$, so we have that $\eta(1,z') \in [z]$ for any $z'
\in [z]$. Fix $\varepsilon > 0$, take $z' \in [z]$ such that
\begin{equation} \label{4.91}
c \le \max \widetilde{J}_s(|z'|) < c + \varepsilon^2,
\end{equation}
and define a continuous function $I : {\cal L} \to [c,c +
\varepsilon^2)$ by
\begin{equation} \label{4.92}
I(\eta) = \max \widetilde{J}_s(|\eta(1,z')|).
\end{equation}
Let $\etabar$ be the identity in $\cal L$ (i.e., $\etabar(t,u) =
u$ for all $(t,u) \in [0,1] \times S$), and note that
\begin{equation} \label{4.93}
I(\etabar) < c + \varepsilon^2 \le \inf_{\cal L} I +
\varepsilon^2.
\end{equation}
Applying the Ekeland's principle, we get an $\eta_0 \in {\cal L}$
such that
\begin{gather}
\label{4.94} I(\eta_0) \le I(\etabar),\\[5pt] \label{4.95}
\delta(\eta_0,\etabar) \le \varepsilon,\\[5pt] \label{4.96}
I(\eta) \ge I(\eta_0) - \varepsilon\, \delta(\eta,\eta_0) \quad
\forall \eta \in {\cal L}.
\end{gather}
Let $C = \set{u \in |\eta_0(1,z')| : \widetilde{J}_s(u) =
I(\eta_0)}$. Since $\widetilde{J}_s$ satisfies (PS) and $c \le
\widetilde{J}_s(u) < c + \varepsilon^2$ for all $u \in C$, it is
enough to show that there is a $u_\varepsilon \in C$ such that
$\norm{\widetilde{J}_s'(u_\varepsilon)}{} \le 4 \varepsilon$.
Indeed, in view of \eqref{4.95}, any such point necessarily
satisfies $\dist{u_\varepsilon}{{\cal O}} \le
\dist{u_\varepsilon}{|z'|} \le \varepsilon$.

Suppose now that $\norm{\widetilde{J}_s'(u)}{} > 4 \varepsilon$
for all $u \in C$. Applying Lemma 3.7 of Ghoussoub \cite{Gh2}, we
get $t_0 > 0$, $\alpha \in C([0,t_0) \times S; S)$, and $g \in
C(S; [0,1])$ such that
\begin{enumeratepr}
\item \label{4.961} $\rho(\alpha(t,u),u) \le \dfrac{3t}{2} \quad
\forall (t,u) \in [0,t_0) \times S$,
\item \label{4.962} $\widetilde{J}_s(\alpha(t,u)) - \widetilde{J}_s(u)
\le - 2 \varepsilon\, g(u)\, t \quad \forall (t,u) \in [0,t_0)
\times S$,
\item \label{4.963} $g(u) = 1 \quad \forall u \in C$.
\end{enumeratepr}
For $0 < \lambda < t_0$, let $\eta_\lambda(t,u) = \alpha(t
\lambda,\eta_0(t,u))$. Clearly, $\eta_\lambda \in {\cal L}$. Since
$\delta(\eta_\lambda,\eta_0) \le \dfrac{3t \lambda}{2} \le
\dfrac{3 \lambda}{2}$ by \eqref{4.961} above, \eqref{4.96} gives
$I(\eta_\lambda) \ge I(\eta_0) - \dfrac{3 \varepsilon
\lambda}{2}$. Since $|z'|$ is compact, there is a $u_\lambda \in
|z'|$ such that $\widetilde{J}_s(\eta_\lambda(1,u_\lambda)) =
I(\eta_\lambda)$, so we have
\begin{equation} \label{4.97}
\widetilde{J}_s(\eta_\lambda(1,u_\lambda)) -
\widetilde{J}_s(\eta_0(1,u)) \ge - \frac{3 \varepsilon \lambda}{2}
\quad \forall u \in |z'|.
\end{equation}
In particular, if $u_0$ is any cluster point of $\seq{u_\lambda}$
as $\lambda \to 0$, then $\eta_0(1,u_0) \in C$, and hence
\begin{equation} \label{4.971}
g(\eta_0(1,u_0)) = 1
\end{equation}
by \eqref{4.963}. On the other hand,
\begin{equation} \label{4.98}
\begin{split}
\widetilde{J}_s(\eta_\lambda(1,u_\lambda)) -
\widetilde{J}_s(\eta_0(1,u_\lambda)) & =
\widetilde{J}_s(\alpha(\lambda,\eta_0(1,u_\lambda)) -
\widetilde{J}_s(\eta_0(1,u_\lambda))\\[5pt] & \le - 2 \varepsilon
\lambda\, g(\eta_0(1,u_\lambda))
\end{split}
\end{equation}
by \eqref{4.962}, and combining this with \eqref{4.97} gives
\begin{equation} \label{4.991}
g(\eta_0(1,u_\lambda)) \le \frac{3}{4},
\end{equation}
which contradicts \eqref{4.971}.
\end{proof}

Note that according to definition \eqref{4.2} there is no path in
${\cal O}_b$ joining $\varphi_1$ and $- \varphi_1$ if $b \le c(s)$.
Conversely, we have

\begin{lemma} \label{T11}
If $b > c(s)$, then ${\cal O}_b$ is path-connected.
\end{lemma}

\begin{proof}
We will show that every point $u \in {\cal O}_b$ can be connected
to $\varphi_1$ by a path in ${\cal O}_b$. Denote by ${\cal O}_u$
the component of ${\cal O}_b$ containing $u$. By Lemma 3.6 of
\cite{CudeFiGo}, $d = \inf \widetilde{J}_s(\overline{{\cal O}_u})$
is achieved at a critical point $u_0 \in {\cal O}_u$. Since
${\cal O}_u$ is path-connected by Lemma 3.5 of \cite{CudeFiGo},
it is enough to show that $u_0$ can be connected to $\varphi_1$
by a path in ${\cal O}_b$.

If $u_0 \le 0$, then $u_0 = - \varphi_1$ since every critical
point of $\widetilde{J}_s$ other than $\pm \varphi_1$ changes
sign, and $- \varphi_1$ can be connected to $\varphi_1$ by a
path in ${\cal O}_b$ since $b > c(s)$. So suppose that $u_0^+ \ne 0$,
and let
\begin{equation} \label{35}
u_t = \frac{u_0^+ - (1-t)\, u_0^-}{\norm{u_0^+ - (1-t)\, u_0^-}{p}},
\quad t \in [0,1].
\end{equation}
Taking $v = u_0^\pm$ in
\begin{equation} \label{36}
\int_\Omega |\nabla u_0|^{p-2}\, \nabla u_0 \cdot \nabla v - s\,
(u_0^+)^{p-1}\, v = d \int_\Omega |u_0|^{p-2}\, u_0\, v
\end{equation}
gives
\begin{equation} \label{37}
\int_\Omega |\nabla u_0^+|^p - s\, (u_0^+)^p = d \int_\Omega
(u_0^+)^p, \qquad \int_\Omega |\nabla u_0^-|^p = d \int_\Omega
(u_0^-)^p,
\end{equation}
so
\begin{equation} \label{38}
\widetilde{J}_s(u_t) = \frac{\displaystyle{\int_\Omega |\nabla
u_0^+|^p + (1-t)^p\, |\nabla u_0^-|^p - s\,
(u_0^+)^p}}{\displaystyle{\int_\Omega (u_0^+)^p + (1-t)^p\,
(u_0^-)^p}} = d.
\end{equation}
Thus $u_1 = \dfrac{u_0^+}{\norm{u_0^+}{p}} \in {\cal O}_u$
achieves the infimum and hence is a critical point of
$\widetilde{J}_s$. Since $u_1 \ge 0$, $u_1 = \varphi_1$, so
$u_t$ is a path in ${\cal O}_b$ connecting $u_0$ to $\varphi_1$.
\end{proof}

We are now ready to give the

\begin{proof}[Proof of Theorem \ref{T2}]
\eqref{10}. Since
\begin{equation} \label{14}
\lambda_1 = \min_{u \in S} \int_\Omega |\nabla u|^p,
\end{equation}
\eqref{4.1} gives
\begin{equation} \label{15}
\widetilde{J}_{a-b}(u) \ge \lambda_1 - \max \set{a,b} + b > b
\quad \forall u \in S,
\end{equation}
so $\widetilde{J}_{a-b}^b = \emptyset$ and the conclusion follows
from Lemma \ref{T1}.

\eqref{11}. By Lemma \ref{T1},
\begin{equation} \label{24}
C_q(I_{(a,b)},0) \cong \widetilde{H}_{q-1}(\widetilde{J}_{a-b}^b).
\end{equation}
We apply Corollary \ref{T8} with $\beta = b$ and ${\cal O} =
\widetilde{J}_{a-b}^b$. The only critical point of
$\widetilde{J}_{a-b}$ in $\widetilde{J}_{a-b}^b$ is the global
minimizer $\varphi_1$, and
\begin{equation} \label{18}
\widetilde{J}_{a-b}(u) > \widetilde{J}_{a-b}(\varphi_1) \quad
\forall u \in X \setminus \set{\varphi_1}
\end{equation}
by Lemma \ref{T4}. Since $\widetilde{J}_{a-b}(\varphi_1) =
\lambda_1 - (a - b) < b$, there is an $\varepsilon > 0$ such that
\begin{equation} \label{17}
\widetilde{J}_{a-b}(u) \le b \quad \forall u \in
\widetilde{B}_\varepsilon(\varphi_1)
\end{equation}
by the continuity of $\widetilde{J}_{a-b}$, and an argument
similar to that in the proof of Lemma 2.9 of \cite{CudeFiGo}
shows that
\begin{equation} \label{19}
\inf_{u \in \bdry{\widetilde{B}_\varepsilon(\varphi_1)}}
\widetilde{J}_{a-b}(u) > \widetilde{J}_{a-b}(\varphi_1).
\end{equation}
Thus
\begin{equation} \label{19.1}
\widetilde{H}_{q-1}(\widetilde{J}_{a-b}^b) = 0.
\end{equation}

\eqref{12}. By Proposition 2.3 and Lemma 2.9 of \cite{CudeFiGo},
there are $\varepsilon > 0$ and $\beta \in (\lambda_1,b]$ such
that
\begin{gather}
\label{25} \widetilde{J}_{a-b}(u) > \widetilde{J}_{a-b}(-
\varphi_1) = \lambda_1 \quad \forall u \in \widetilde{B}_\varepsilon(-
\varphi_1) \setminus \set{- \varphi_1},\\[5pt]
\label{25.1} \inf_{u \in \bdry{\widetilde{B}_\varepsilon(-
\varphi_1)}} \widetilde{J}_{a-b}(u) > \beta.
\end{gather}
Since $\widetilde{J}_{a-b}$ has no critical values in $[\beta,b]$,
\begin{equation} \label{26}
\widetilde{H}_{q-1}(\widetilde{J}_{a-b}^b) \cong
\widetilde{H}_{q-1}(\widetilde{J}_{a-b}^\beta)
\end{equation}
by Lemma \ref{T5}. As \eqref{25.1} implies that ${\cal O}_1 =
\widetilde{J}_{a-b}^\beta \setminus \widetilde{B}_\varepsilon(-
\varphi_1)$ and ${\cal O}_2 = \widetilde{J}_{a-b}^\beta \setminus
{\cal O}_1 = \widetilde{J}_{a-b}^\beta \cap
\widetilde{B}_\varepsilon(- \varphi_1)$ are isolated in
$\widetilde{J}_{a-b}^\beta$,
\begin{equation} \label{27}
\widetilde{H}_{q-1}(\widetilde{J}_{a-b}^\beta) \cong
\widetilde{H}_{q-1}({\cal O}_1) \oplus H_{q-1}({\cal O}_2),
\end{equation}
and
\begin{equation} \label{28}
\widetilde{H}_{q-1}({\cal O}_1) = 0, \qquad H_{q-1}({\cal O}_2) =
\delta_{q1}\, \Z \oplus \widetilde{H}_{q-1}({\cal O}_2) =
\delta_{q1}\, \Z
\end{equation}
by Corollary \ref{T8}.

\eqref{13}. By Lemma \ref{T1},
\begin{equation} \label{34}
C_q(I_{(a,b)},0) \cong \widetilde{H}_{q-1}({\cal O}_b),
\end{equation}
and the conclusion follows since ${\cal O}_b$ is path-connected by
Lemma \ref{T11}.
\end{proof}

\begin{proof}[Proof of Proposition \ref{T6}]
If both $(a_0,b_0)$ and $(a_1,b_1)$ are below $\Cl{1}$, then the
conclusion follows from part \eqref{10} of Proposition \ref{T2}.
So, by Lemma \ref{T1}, it suffices to show that for any $(a_0,b_0)
\notin \Sigma_p$ that lies above $\Cl{1}$,
\begin{equation} \label{29}
\widetilde{H}_\ast(\widetilde{J}_{a-b}^b) \cong
\widetilde{H}_\ast(\widetilde{J}_{a_0-b_0}^{b_0})
\end{equation}
for all $(a,b)$ sufficiently close to $(a_0,b_0)$.

Choose $\varepsilon > 0$ so small that $B_{4 \varepsilon}(a_0,b_0)
\cap \Sigma_p = \emptyset$, which is possible since $\Sigma_p$ is
closed, and suppose that $|a - a_0| + |b - b_0| \le \varepsilon$.
Then
\begin{equation} \label{30}
\left|\widetilde{J}_{a-b}(u) - \widetilde{J}_{a_0-b_0}(u)\right| =
|(a - a_0) - (b - b_0)| \int_\Omega (u^+)^p \le \varepsilon \quad
\forall u \in S,
\end{equation}
so we have the inclusions
\begin{equation} \label{31}
\xymatrix{{\widetilde{J}}_{a_0-b_0}^{b - \varepsilon}
\ar@{^{(}->}[r]^(0.54){i_1} & {\widetilde{J}}_{a-b}^b
\ar@{^{(}->}[r]^(0.46){i_2} & {\widetilde{J}}_{a_0-b_0}^{b +
\varepsilon} \ar@{^{(}->}[r]^{i_3} & {\widetilde{J}}_{a-b}^{b + 2
\varepsilon},}
\end{equation}
which induce homomorphisms
\begin{equation} \label{32}
\xymatrix{{\widetilde{H}}_\ast(\widetilde{J}_{a_0-b_0}^{b -
\varepsilon}) \ar[r]^(0.54){i_{1 \ast}} &
{\widetilde{H}}_\ast(\widetilde{J}_{a-b}^b) \ar[r]^(0.49){i_{2
\ast}} & {\widetilde{H}}_\ast(\widetilde{J}_{a_0-b_0}^{b +
\varepsilon}) \ar[r]^(0.51){i_{3 \ast}} &
{\widetilde{H}}_\ast(\widetilde{J}_{a-b}^{b + 2 \varepsilon}).}
\end{equation}
Since the points $(a_0 - b_0 + b - \varepsilon,b - \varepsilon)$
and $(a_0 - b_0 + b + \varepsilon,b + \varepsilon)$ (resp. $(a,b)$
and $(a + 2 \varepsilon,b + 2 \varepsilon)$) are in $B_{4
\varepsilon}(a_0,b_0)$, $\widetilde{J}_{a_0-b_0}$ (resp.
$\widetilde{J}_{a-b}$) has no critical values in $[b -
\varepsilon,b + \varepsilon]$ (resp. $[b,b + 2 \varepsilon]$), so
$i_{2 \ast}\, i_{1 \ast}$ (resp. $i_{3 \ast}\, i_{2 \ast}$) is an
isomorphism by Lemma \ref{T5}. Thus $i_{2 \ast}$ is an
isomorphism. Finally,
\begin{equation} \label{33}
\widetilde{H}_\ast(\widetilde{J}_{a_0-b_0}^{b + \varepsilon})
\cong \widetilde{H}_\ast(\widetilde{J}_{a_0-b_0}^{b_0})
\end{equation}
as $\widetilde{J}_{a_0-b_0}$ has no critical values between
$b_0$ and $b + \varepsilon$.
\end{proof}

\section{Proof of Theorem \ref{T9}} \label{S3}

Let
\begin{gather}
\label{54} f_\pm(x,t) = \begin{cases}
f(x,t) & \text{if } \pm t \ge 0,\\[5pt]
0 & \text{otherwise},
\end{cases} \qquad F_\pm(x,t) = \int_0^t f_\pm(x,s)\, ds,\\[5pt]
\label{55} \Phi_\pm(u) = \int_\Omega |\nabla u|^p - p\, F_\pm(x,u).
\end{gather}
If $u$ is a critical point of $\Phi_+$, taking $v = u^-$ in
\begin{equation} \label{56}
\ip{\Phi_+'(u)}{v} = \int_\Omega |\nabla u|^{p-2}\, \nabla u
\cdot \nabla v - f_+(x,u)\, v = 0
\end{equation}
shows that $u^- = 0$, so $u = u^+$ is also a critical point of
$\Phi$ with critical value $\Phi(u) = \Phi_+(u)$. Furthermore,
$u \in L^\infty(\Omega) \cap C^1(\Omega)$ by Anane \cite{An} and
di Benedetto \cite{diBe}, so it follows from the Harnack
inequality (Theorem 1.1 of Trudinger \cite{Tr}) that either $u >
0$ or $u \equiv 0$. Similarly, nontrivial critical points of
$\Phi_-$ are negative solutions of \eqref{39}.

The following lemma is known, but for completeness we give a
simple proof.

\begin{lemma} \label{T13}
The mapping $u \mapsto |u|$ is continuous on $X$. Hence the
mappings $u \mapsto u^\pm = \dhalf\, \big(|u| \pm u\big)$ are also
continuous on $X$.
\end{lemma}

\begin{proof}
Since $\nabla |u| = \sgn u\, \nabla u$ a.e., $\norm{|u|}{} =
\norm{u}{}$, so by uniform convexity it suffices to prove norm to
weak continuity, i.e., $u_j \to u$ implies $|u_j| \rightharpoonup
|u|$. We have $|u_j| \to |u|$ in $L^p(\Omega)$, so by weak
compactness $|u_j| \rightharpoonup z$ in X for a subsequence.
Hence $z =|u|$ and $|u_j| \rightharpoonup |u|$ for the whole
sequence.
\end{proof}

\begin{lemma} \label{T12}
Nontrivial local minimizers of $\Phi_\pm$ are also local
minimizers of $\Phi$.
\end{lemma}

\begin{proof}
We only consider a local minimizer $u_0 > 0$ of $\Phi_+$ as the
argument for $\Phi_-$ is similar. We have to show that for every
sequence $u_j \to u_0$ in $X$, $\Phi(u_j) \ge \Phi(u_0)$ for
sufficiently large $j$. By \eqref{40},
\begin{equation} \label{70}
|F(x,t)| \le C\, |t|^p \quad \forall t
\end{equation}
for some constant $C > 0$, so
\begin{align} \label{69}
\Phi(u_j) & = \int_\Omega |\nabla u_j^+|^p - p\, F(x,u_j^+) +
\int_\Omega |\nabla u_j^-|^p - p\, F(x,u_j^-)\\[5pt]
& \ge \Phi_+(u_j^+) + \norm{u_j^-}{}^p - C\, \norm{u_j^-}{p}^p.
\end{align}
Since $u_j^+ \to u_0^+ = u_0$ by Lemma \ref{T13},
\begin{equation} \label{71}
\Phi_+(u_j^+) \ge \Phi_+(u_0) = \Phi(u_0)
\end{equation}
for sufficiently large $j$. We will show that
\begin{equation} \label{67}
\norm{u_j^-}{}^p \ge C\, \norm{u_j^-}{p}^p, \quad j \text{ large}.
\end{equation}

First we note that the measure of the set $\Omega_j = \set{x \in
\Omega : u_j(x) < 0}$ goes to zero. To see this, given $\varepsilon
> 0$, take a compact subset $\Omega^\varepsilon$ of $\Omega$ such
that $|\Omega
\setminus \Omega^\varepsilon| < \varepsilon$ and let
$\Omega^\varepsilon_j = \Omega^\varepsilon \cap \Omega_j$. Then
\begin{equation} \label{64}
\norm{u_j - u_0}{p}^p \ge \int_{\Omega^\varepsilon_j} |u_j - u_0|^p
\ge \int_{\Omega^\varepsilon_j} u_0^p \ge c^p\,
|\Omega^\varepsilon_j|
\end{equation}
where $c = \min_{\Omega^\varepsilon} u_0 > 0$, so
$|\Omega^\varepsilon_j| \to 0$. Since $\Omega_j \subset
\Omega^\varepsilon_j \cup (\Omega \setminus \Omega^\varepsilon)$ and
$\varepsilon > 0$ is arbitrary, the claim follows.

If \eqref{67} does not hold, setting $\tilde{u}_j =
\dfrac{u_j^-}{\norm{u_j^-}{p}}$, $\norm{\tilde{u}_j}{}$ is bounded
for a subsequence, so $\tilde{u}_j \to \tilde{u}$ in $L^p(\Omega)$
and a.e. for a further subsequence, where $\norm{\tilde{u}}{p} =
1$ and $\tilde{u} \ge 0$. But then $\Omega_\mu = \set{x \in \Omega
: \tilde{u}(x) \ge \mu}$ has positive measure for all sufficiently
small $\mu > 0$ and
\begin{equation} \label{68}
\norm{\tilde{u}_j - \tilde{u}}{p}^p \ge \int_{\Omega_\mu \setminus
\Omega_j} |\tilde{u}_j - \tilde{u}|^p = \int_{\Omega_\mu \setminus
\Omega_j} \tilde{u}^p \ge \mu^p\, (|\Omega_\mu| - |\Omega_j|),
\end{equation}
a contradiction.
\end{proof}

\begin{lemma} \label{T10}
If $(a_0,b_0),\, (a,b) \notin \Sigma_p$, then there are $R > \rho
> 0$ and $\widetilde{\Phi} \in C^1(X,\R)$ such that
\begin{equation} \label{53}
\widetilde{\Phi} = \begin{cases}
I_{(a_0,b_0)} & \text{in } \overline{B}_{\rho/2},\\[5pt]
\Phi & \text{in } \overline{B}_R \setminus B_\rho,\\[5pt]
I_{(a,b)} & \text{in } X \setminus B_{2R}
\end{cases}
\end{equation}
and $0$ is the only critical point of $\Phi$ and
$\widetilde{\Phi}$ in $\overline{B}_\rho \cup (X \setminus B_R)$,
where $B_\rho = \bigset{u \in X : \norm{u}{} < \rho}$. In
particular,
\begin{equation} \label{53.1}
C_\ast(\widetilde{\Phi},0) = C_\ast(I_{(a_0,b_0)},0), \qquad
C_\ast(\widetilde{\Phi},\infty) = C_\ast(I_{(a,b)},0).
\end{equation}
\end{lemma}

\begin{proof}
Let $g_0(x,t) = f(x,t) - a_0\, (t^+)^{p-1} + b_0\,
(t^-)^{p-1},\, g(x,t) = f(x,t) - a\, (t^+)^{p-1} + b\,
(t^-)^{p-1},\, G_0(x,t) = {\displaystyle \int_0^t}
g_0(x,s)\, ds,\, G(x,t) = {\displaystyle \int_0^t} g(x,s)\,
ds$, and
\begin{equation} \label{43}
\Psi_0(u) = - p \int_\Omega G_0(x,u), \qquad \Psi(u) = - p
\int_\Omega G(x,u),
\end{equation}
so that
\begin{equation} \label{44}
\Phi(u) = I_{(a_0,b_0)}(u) + \Psi_0(u) = I_{(a,b)}(u) +
\Psi(u).
\end{equation}
Since $(a_0,b_0),\, (a,b) \notin \Sigma_p$, $I_{(a_0,b_0)}$
and $I_{(a,b)}$ satisfy (PS) and have no critical points on
$S_1 = \bdry{B_1}$, so
\begin{equation} \label{45}
\delta_0 = \inf_{S_1} \norm{I_{(a_0,b_0)}'}{} > 0, \qquad
\delta = \inf_{S_1} \norm{I_{(a,b)}'}{} > 0.
\end{equation}
By homogeneity,
\begin{equation} \label{46}
\inf_{S_\rho} \norm{I_{(a_0,b_0)}'}{} = \rho^{p-1}\,
\delta_0, \qquad \inf_{S_R} \norm{I_{(a,b)}'}{} = R^{p-1}\,
\delta,
\end{equation}
while it follows from \eqref{40} that
\begin{gather}
\label{47} \sup_{S_\rho} |\Psi_0| = o(\rho^p), \qquad
\sup_{S_R} |\Psi| = o(R^p),\\[5pt]
\label{48} \sup_{S_\rho} \norm{\Psi_0'}{} = o(\rho^{p-1}), \qquad
\sup_{S_R} \norm{\Psi'}{} = o(R^{p-1})
\end{gather}
as $\rho \to 0$ and $R \to \infty$, so
\begin{equation} \label{49}
\inf_{S_\rho} \norm{\Phi'}{} \ge \rho^{p-1} (\delta_0 +
o(1)) > 0, \quad \halfquad \inf_{S_R} \norm{\Phi'}{} \ge
R^{p-1} (\delta + o(1)) > 0
\raisetag{22pt}
\end{equation}
for all sufficiently small $\rho > 0$ and sufficiently
large $R > \rho$. Take smooth functions $\varphi_0,\,
\varphi : [0,\infty) \to [0,1]$ such that
\begin{equation} \label{50}
\varphi_0(t) = \begin{cases}
1 & \text{if } 0 \le t \le 1/2,\\[5pt]
0 & \text{if } t \ge 1,
\end{cases} \qquad \varphi(t) = \begin{cases}
0 & \text{if } t \le 1,\\[5pt]
1 & \text{if } t \ge 2,
\end{cases}
\end{equation}
and set
\begin{equation} \label{51}
\widetilde{\Phi}(u) = \Phi(u)
- \varphi_0(\norm{u}{}/\rho)\, \Psi_0(u)
- \varphi(\norm{u}{}/R)\, \Psi(u).
\end{equation}
Since
\begin{equation} \label{52}
\norm{d\left(\varphi_0(\norm{u}{}/\rho)\right)}{} =
O(\rho^{-1}), \qquad
\norm{d\left(\varphi(\norm{u}{}/R)\right)}{} = O(R^{-1}),
\end{equation}
\eqref{49} holds with $\Phi$ replaced by $\widetilde{\Phi}$
also. The conclusion follows.
\end{proof}

We are now ready to prove Theorem \ref{T9}. Denote by
$\widetilde{\Phi}_\pm$ the modified functionals obtained by
applying Lemma \ref{T10}.

\eqref{60}. If $a_0 < \lambda_1 < a$, then $(a_0,0)$ is below
$\Cl{1}$ and $(a,0)$ is between $\Cl{1}$ and $\Cu{1}$, so
\begin{gather}
\label{57} C_q(\widetilde{\Phi}_+,0) = C_q(I_{(a_0,0)},0) =
\delta_{q0}\, \Z,\\[5pt]
\label{57.1} C_q(\widetilde{\Phi}_+,\infty) = C_q(I_{(a,0)},0)
= 0 \quad \forall q
\end{gather}
by Lemma \ref{T10} and Theorem \ref{T2}. Similarly, if $a_0 >
\lambda_1 > a$, then
\begin{equation} \label{58}
C_q(\widetilde{\Phi}_+,0) = 0 \quad \forall q, \qquad
C_q(\widetilde{\Phi}_+,\infty) = \delta_{q0}\, \Z.
\end{equation}
In either case $C_0(\widetilde{\Phi}_+,0) \not \cong
C_0(\widetilde{\Phi}_+,\infty)$, so $\widetilde{\Phi}_+$ must
have a nontrivial critical point.

\eqref{61}. Follows similarly by working with
$\widetilde{\Phi}_-$.

\eqref{59}. In view of \eqref{60} and \eqref{61}, it only
remains to consider the case where both points are above $\Cu{1}$
and on opposite sides of $\C{2}$. If $(a_0,b_0)$ is between
$\Cu{1}$ and $\C{2}$ and $(a,b)$ is above $\C{2}$, then
\begin{equation} \label{65}
C_q(\widetilde{\Phi},0) = \delta_{q1}\, \Z, \qquad
C_1(\widetilde{\Phi},\infty) = 0.
\end{equation}
The other case is similar.

\eqref{62}. Follows from \eqref{60} and \eqref{61}.

\eqref{63}. Since $a,\, b < \lambda_1$, $\Phi_\pm$ are bounded
from below and coercive, and hence have global minimizers
$u_0^\pm$. As $a_0,\, b_0 > \lambda_1$,
\begin{equation} \label{72}
\Phi_\pm(u_0^\pm) \le \Phi_\pm(\pm\, \rho\halfthinspace \phi_1) <
0 = \Phi_\pm(0), \quad \rho > 0 \text{ small},
\end{equation}
so $u_0^\pm \ne 0$. Thus $u_0^\pm$ are local minimizers of $\Phi$,
and hence also of $\widetilde{\Phi}$, by Lemma \ref{T12}. Now a
standard mountain-pass argument yields a critical point $u_1 \ne
u_0^\pm$ of $\widetilde{\Phi}$ with
\begin{equation} \label{73}
C_1(\widetilde{\Phi},u_1) \ne 0,
\end{equation}
which is nontrivial since
\begin{equation} \label{74}
C_1(\widetilde{\Phi},0) = 0. \qed
\end{equation}

\end{document}